\documentclass[twoside]{article}
\usepackage{graphics, graphicx, latexsym, amssymb, amsmath}
\newtheorem{theo}{Theorem}
\newtheorem{coro}{Corollary}
\begin{document}
\title{The area generating function \\ for simple-2-column polyominoes \\ with hexagonal cells}

\author{Svjetlan Fereti\'{c} \footnote{e-mail: svjetlan.feretic@gradri.hr} \\ 
Faculty of Civil Engineering, University of Rijeka, \\ 
Viktora Cara Emina 5, HR-51\,000 Rijeka, Croatia
\and
Nenad Trinajsti\'{c} \footnote{e-mail: trina@irb.hr} \\
The Rugjer Bo\v{s}kovi\'{c} Institute,\\
P. O. Box 180, HR-10\,002 Zagreb, Croatia}

\maketitle

\begin{abstract}
Column-convex polygons were first counted by area several decades ago, and the result was found to be a simple, rational, generating function. In this chapter we generalize that result. Let a $p$-column polyomino be a polyomino whose columns can have $1,\ 2,\ldots,\ p$ connected components. Then column-convex polygons are equivalent to 1-convex polyominoes. The area generating function of even the simplest generalization, namely to 2-column polyominoes, is unlikely to be solvable. We therefore define a class of polyominoes which interpolates between column-convex polygons and 2-column polyominoes. We derive the area generating function of that  class, using an extension of an existing algorithm. The growth constant of the new class is greater than the growth constant of column-convex polyominoes. A rather tight lower bound on the growth constant complements a compelling numerical analysis.
\end{abstract}

\section{Introduction}

The enumeration of polyominoes is a topic of great interest to chemists, physicists and combinatorialists alike \cite{book}. In chemical terms, any polyomino (with hexagonal cells and with no internal holes) is a possible benzenoid hydrocarbon. In combinatorics, polyominoes are of interest in their own right because several polyomino models have mathematically appealing exact solutions. Furthermore, they are also relevant to various problems of tilings \cite{GS87}. Polyominoes are extensively studied, in one form or another, in a wide variety of problems of interest to physicists. In particular, we note their investigation under the name {\em lattice animals}, in the study of percolation \cite{Ha05, J01}, in the graphical representation of the Ising model, and its extension to the Potts model, and in the study of the properties of branched polymers \cite{HN83, PS89, V85}. 

They are also a representative of a class of problems that appear to be unsolvable---notably the enumeration, by area or perimeter, of self-avoiding polygons, polyominoes and other classes of graphs that are relevant to the Ising model and to percolation.  The principal line of attack on such problems is to simplify them until they are solvable, in the hope that the essential physics is not destroyed in the process. That is the approach taken in this chapter, where the model proposed, while still solvable, is closer to the ultimate problem of full polyomino enumeration than has previously been attained. Further, by use of 
Bousquet-M\'{e}lou's \cite{Bousquet} and Svrtan's~\cite{Svrtan} upgraded version of the Temperley methodology \cite{Temperley}, we give the solution of one problem previously out of reach due to its complexity. This development may spur further advances in this class of problem.

One frequently cited polyomino model is that of column-convex polygons~\footnote{We distinguish between polygons and polyominoes in that the former cannot have internal holes. As a consequence, the {\it perimeter} generating function for polygons has a non-zero radius of convergence, whereas for polyominoes the radius of convergence is zero.}. There exist two main versions of column-convex polygons: the first composed of square cells and the second of hexagonal cells. Both versions have a rational area generating function. For the version with square cells, the area generating function was found independently by P\'{o}lya \cite{Polya} in 1938 \textit{or} 1969 and by Temperley \cite{Temperley} in 1956. That was perhaps the earliest major result in polyomino enumeration. For the version with hexagonal cells, the area generating function was found by Klarner in 1967 \cite{Klarner}. The growth constant of square-celled column-convex polygons is $\mu=3.205569\ldots \:$, while the growth constant of hexagonal-celled column-convex polygons is $\mu=3.863130\ldots\ $. (By the growth constant we mean the limit $\lim_{n\rightarrow\infty} \sqrt[n]{a_n}$, where $a_n$ denotes the number of $n$-celled elements in a given set of polyominoes.) In both cases the area generating function is a simple pole, so that $a_n \sim const. \times \mu^n.$

There exist some models which are supersets of column-convex polygons and are still solvable. Those models are called $m$\textit{-convex polygons} \cite{m-convex}, \textit{prudent polygons} \cite{DG08}, \textit{cheesy polyominoes} \cite{semi}, \textit{polyominoes with cheesy blocks} \cite{semiblo}, \textit{column-subconvex polyominoes} \cite{undi}, and \textit{simple-2-column polyominoes} \cite{simplex-duplex}. The former two models can be enumerated by perimeter and area, whereas the latter four models have been enumerated only by area. 

This chapter is a kind of companion to \cite{simplex-duplex}. Namely, the enumeration of simple-2-column polyominoes by area was first done in \cite{simplex-duplex}. In the present chapter, we again enumerate simple-2-column polyominoes by area, but this time the cells are hexagons, whereas in \cite{simplex-duplex} the cells were squares. As mentioned above, we make use of Bousquet-M\'{e}lou's \cite{Bousquet} and Svrtan's~\cite{Svrtan} upgraded version of Temperley's methodology \cite{Temperley}. The computations are rather long and intricate. Also, these computations are very similar to those which are done detailedly in \cite{simplex-duplex}. Therefore in this chapter we only give an outline of the proof, though with enough detail that the method may be applied by others to new problems. Incidentally, in \cite{simplex-duplex}, simple-2-column polyominoes are called by their original name (coined by S. Fereti\'{c}): \textit{simplex-duplex polyominoes}. The name simple-2-column polyominoes was suggested subsequently (during the preparation of \cite{FG09}) by Tony Guttmann. 

In Section 2, we discuss the chemical relevance of polyominoes with hexagonal cells. In Section 3, we define the model. In Section 4 we give the formula for $G(q,w)$, a generating function for simple-2-column polyominoes, in which the variable $q$ is conjugate to the area and $w$ is conjugate to the number of two-component columns of the polyomino. A truncated version of the proof is given in Section 5. In Section 6, we discuss the asymptotic behaviour of $G(q,w)$, and give a tight lower bound on the growth constant, as well as a very accurate estimate. In Section 7 we conclude, outlining further work prompted by our results. 

Note that our solution essentially gives detailed information only about the area generating function. The additional parameter $w$ counts columns of a certain type. Thus we cannot give perimeter-area phase diagrams which are relevant to the description of vesicle collapse. Indeed, as is shown in \cite{FG09}, the perimeter generating function has zero radius of convergence (as is the usual case for polyominoes), which precludes such a phase diagram in its usual form.

\section{Chemical relevance of polyominoes\\ with hexagonal cells}

Polyominoes with hexagonal cells can be used to model benzenoid hydrocarbons. Benzenoid hydrocarbons are a class of versatile conjugated organic molecules that have constantly been studied by a variety of researchers. This is so because benzenoid hydrocarbons found uses in experimental and theoretical research in chemistry \cite{Clar, Dias, Gutman}, environmental chemistry \cite{Futoma, Neff}, chemical technology \cite{Balaban}, cancer research \cite{Daudel, Gelboin, Harvey}, ecology \cite{Butler, Lunde}, combinatorial chemistry \cite{Knop_12}, graph theory \cite{Harary, Trinajstic_14}, theory of aromaticity \cite{Minkin}, computational chemistry \cite{Balasubramanian, Knop_18, Trinajstic_17, Trinajstic_19}, polymer science \cite{PS89}, \textit{etc.}

Additionally, polyominoes and benzenoid structures under various disguises have been discussed regularly at the Dubrovnik MATH/CHEM/COMP Symposia. For example, this class of structures has been discussed there as fractal benzenoids \cite{Plavsic}, as column-convex animals \cite{Delest, rani, everybody}, in algebraic studies of benzenoids \cite{Cyvin}, in leapfrog transformation \cite{Babic}, maximum matchings and eigenvalues of benzenoid graphs \cite{Fajtlowicz}, in phototoxicity \cite{Estrada}, ring-currents in benzenoid hydrocarbons \cite{Mallion}, \textit{etc.}

Apparently research on various aspects of this class of structures continues to date and the present chapter is in this direction. However, before we proceed, we wish to point to perhaps the first scientific study on hexagonal structures. This study is due to Croatian polymath Rugjer Bo\v{s}kovi\'{c} (1711--1787). He studied the structure of honeycomb. Honeycomb is made of hexagonal cells. Bo\v{s}kovi\'{c} hypothesized that the hexagonal lattice is not chosen by chance. By using hexagonal lattice, bees economize on building material. Bo\v{s}kovi\'{c} published this study, entitled \textit{De apium cellulis}, in 1760.

\section{Definition of the model}

There are three regular tilings of the Euclidean plane, namely the triangular tiling, the square tiling, and the hexagonal tiling. We adopt the convention that every square or hexagonal tile has two horizontal edges. In a regular tiling, a tile is often referred to as a \textit{cell}. A plane figure $P$ is a \textit{polyomino} if $P$ is a union of finitely many cells and the interior of $P$ is connected. 
Observe that, if a union of \textit{hexagonal} cells is connected, then it possesses a connected interior as well, as a connected union of hexagonal tiles must be connected through shared edges. Topologically, a connected union of square cells may be connected only at a shared vertex. Such unions are forbidden by the definition of polyominoes however.

Let $P$ and $Q$ be two polyominoes. We consider $P$ and $Q$ to be equal if and only if there exists a translation $f$ such that $f(P)=Q$.

From now on, we concentrate on the hexagonal tiling. When we write ``a polyomino'', we actually mean ``a hexagonal-celled polyomino''.

Given a polyomino $P$, it is useful to partition the cells of $P$ according to their horizontal projection. Each block of that partition is a \textit{column} of $P$. Note that a column of a polyomino is not necessarily a connected set. On the other hand, it may happen that every column of a polyomino $P$ is a connected set.
In this case, the polyomino $P$ is a \textit{column-convex polygon}. See Figure 1.

\begin{figure}
\begin{center}
\includegraphics[width=47.5mm]{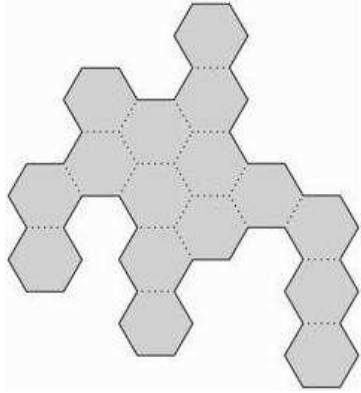}
\caption{A column-convex polygon.}
\end{center}
\end{figure}

By a \textit{2-column polyomino}, we mean a polyomino in which columns with three or more connected components are not allowed. Thus, each column of a 2-column polyomino has either one or two connected components. A \textit{simple-2-column polyomino} is such a 2-column polyomino in which consecutive two-component columns are not allowed. If $c$ is a column of a simple-2-column polyomino, and $c$ is a (left or right) neighbour of a two-component column, then $c$ must be a one-component column. See Figure 2.

\begin{figure}
\begin{center}
\includegraphics[width=70mm]{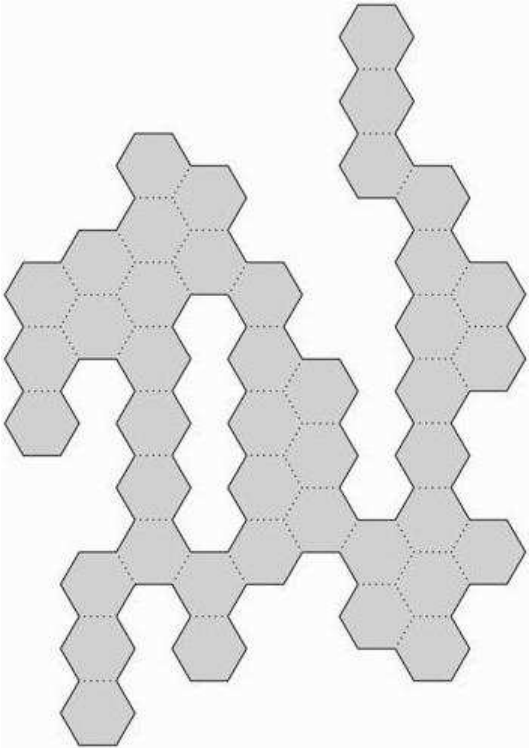} 
\caption{A simple-2-column polyomino.}
\end{center}
\end{figure}

\section{The area generating function for\newline simple-2-column polyominoes\newline with hexagonal cells}

If a polyomino $P$ is made up of $n$ cells, we say that the \textit{area} of $P$ is $n$. Let $ {\mathcal R}$ denote the set of all simple-2-column polyominoes with hexagonal cells.

In Theorem 2 below, we state a formula for the generating function

\begin{displaymath}
G(q,w)=\sum_{{P \in \mathcal{R}}} q^{area\ of\ P} \cdot w^{number\ of\ two-component\ columns\ of\ P}.
\end{displaymath}

\begin{theo} The generating function $G(q,w)$ is given by

\begin{equation}
\label{flag}
G(q,w)=\frac{NUM}{DEN},
\end{equation}

\noindent where

\begin{eqnarray*}
NUM & = & (1-q)^4({\alpha}+{\gamma}+2{\alpha}{\eta}-2{\gamma}{\epsilon})+q^2w(1-q)^2({\iota}+{\lambda}-{\alpha}{\kappa}-{\alpha}{\mu} \\
& & \mbox{}+{\beta}{\iota}+{\beta}{\lambda}-{\gamma}{\kappa}-{\gamma}{\mu}
+{\delta}{\iota}+{\delta}{\lambda}-2{\epsilon}{\lambda}+2{\eta}{\iota} +2{\alpha}{\zeta}{\lambda}-2{\alpha}{\eta}{\kappa} \\
& & \mbox{}-2{\alpha}{\eta}{\mu}+2{\alpha}{\theta}{\lambda}-2{\beta}{\epsilon}{\lambda}
+2{\beta}{\eta}{\iota}+2{\gamma}{\epsilon}{\kappa}+2{\gamma}{\epsilon}{\mu}
-2{\gamma}{\zeta}{\iota}-2{\gamma}{\theta}{\iota} \\
& & \mbox{}-2{\delta}{\epsilon}{\lambda}+2{\delta}{\eta}{\iota}) +4q^2w(1-q)({\alpha}{\lambda} -{\gamma}{\iota}), \\
DEN & = & (1-q)^4(1-{\beta}+{\delta}-{\epsilon}+{\eta}-{\alpha}{\zeta}+{\alpha}{\theta}
+{\beta}{\epsilon}-{\beta}{\eta}+{\gamma}{\zeta}
-{\gamma}{\theta} \\
& & \mbox{}-{\delta}{\epsilon}+{\delta}{\eta})-2(1-q)^3({\gamma}+{\alpha}{\eta}-{\gamma}{\epsilon}) \\
& & \mbox{}-2q^2w(1-q)^2({\kappa}-{\beta}{\mu}+{\delta}{\kappa}-{\epsilon}{\kappa}+{\zeta}{\iota} 
-{\zeta}{\lambda}+{\eta}{\kappa}-{\alpha}{\zeta}{\mu}+{\alpha}{\theta}{\kappa} \\
& & \mbox{}+{\beta}{\epsilon}{\mu}-{\beta}{\eta}{\mu}-{\beta}{\theta}{\iota}
+{\beta}{\theta}{\lambda}+{\gamma}{\zeta}{\mu}-{\gamma}{\theta}{\kappa}
-{\delta}{\epsilon}{\kappa}+{\delta}{\zeta}{\iota}-{\delta}{\zeta}{\lambda} +{\delta}{\eta}{\kappa}) \\
& & \mbox{}-2q^2w(1-q)({\iota}+{\alpha}{\kappa}-{\alpha}{\mu}-{\beta}{\iota}+2{\beta}{\lambda}-2{\gamma}{\kappa}
+{\delta}{\iota}-{\epsilon}{\lambda}+{\eta}{\iota} \\
& & \mbox{}+{\alpha}{\zeta}{\lambda}-{\alpha}{\eta}{\kappa}-{\alpha}{\eta}{\mu}+{\alpha}{\theta}{\lambda}
-{\beta}{\epsilon}{\lambda}+{\beta}{\eta}{\iota}+{\gamma}{\epsilon}{\kappa}+{\gamma}{\epsilon}{\mu}
-{\gamma}{\zeta}{\iota}-{\gamma}{\theta}{\iota} \\
& & \mbox{}-{\delta}{\epsilon}{\lambda}+{\delta}{\eta}{\iota})-4q^2w({\alpha}{\lambda}-{\gamma}{\iota}),
\end{eqnarray*}

\begin{eqnarray*}
\beta & = & \sum_{i=1}^{\infty} \frac{(-3)^{i-1}q^{i^2+2i-2}w^{i-1}}{(1-q)^{2i-2} \cdot \left[\prod_{k=1}^{i-1}(1-q^k)\right]^4 \cdot 
(1-q^i)^{{\scriptstyle{\overline{2}}}\atop}} \ , \\
\gamma & = & \sum_{i=1}^{\infty} \frac{(-3)^{i-1}q^{i^2+3i}w^i}{(1-q)^{2i} \cdot \left[\prod_{k=1}^{i-1}(1-q^k)\right]^4 \cdot 
(1-q^i)^{{\scriptstyle{\overline{3}}}\atop}} \ , \\
\zeta & = & \sum_{i=1}^{\infty} \frac{(-3)^{i-1}q^{i^2+2i-2}w^{i-1}}{(1-q)^{2i-2} \cdot \left[\prod_{k=1}^{i-1}(1-q^k)\right]^4 \cdot 
(1-q^i)^{{\scriptstyle{\overline{2}}}\atop}} \\ 
& & \cdot \left(2i-1+4 \cdot \sum_{k=1}^{i-1}\frac{q^k}{1-q^k} + \frac{\overline{2}q^i}{1-q^i} \right) , 
\end{eqnarray*}

\begin{eqnarray*}
\eta & = & \sum_{i=1}^{\infty} \frac{(-3)^{i-1}q^{i^2+3i}w^i}{(1-q)^{2i} \cdot \left[\prod_{k=1}^{i-1}(1-q^k)\right]^4 \cdot 
(1-q^i)^{{\scriptstyle{\overline{3}}}\atop}} \\ 
& & \cdot \left(2i+4 \cdot \sum_{k=1}^{i-1}\frac{q^k}{1-q^k} + \frac{\overline{3}q^i}{1-q^i} \right) , \\
\kappa & = & \frac{1}{2} \cdot \sum_{i=1}^{\infty} \frac{(-3)^{i-1}q^{i^2+2i-2}w^{i-1}}{(1-q)^{2i-2} \cdot \left[\prod_{k=1}^{i-1}(1-q^k)\right]^4 \cdot 
(1-q^i)^{{\scriptstyle{\overline{2}}}\atop}} \\ 
& & \cdot \left[\left(2i-1+4 \cdot \sum_{k=1}^{i-1}\frac{q^k}{1-q^k} + \frac{\overline{2}q^i}{1-q^i} \right)^2 \right. \\
& & \left. -2i+1+4 \cdot \sum_{k=1}^{i-1}\frac{q^{2k}}{(1-q^k)^2} + \frac{\overline{2}q^{2i}}{(1-q^i)^2}\right] , \\
\lambda & = & \frac{1}{2} \cdot \sum_{i=1}^{\infty} \frac{(-3)^{i-1}q^{i^2+3i}w^i}{(1-q)^{2i} \cdot \left[\prod_{k=1}^{i-1}(1-q^k)\right]^4 \cdot 
(1-q^i)^{{\scriptstyle{\overline{3}}}\atop}} \\ 
& & \cdot \left[\left(2i+4 \cdot \sum_{k=1}^{i-1}\frac{q^k}{1-q^k} + \frac{\overline{3}q^i}{1-q^i} \right)^2 \right. \\
& & \left. -2i+4 \cdot \sum_{k=1}^{i-1}\frac{q^{2k}}{(1-q^k)^2} + \frac{\overline{3}q^{2i}}{(1-q^i)^2} \right] .
\end{eqnarray*}

\emph{In the above formulae, it will be noticed that (a) some of the numbers have an overline, and (b) no result is given for $\alpha$, $\delta$, $\epsilon$, $\theta$, $\iota$, and $\mu$.}

\emph{This is both to save space, and to highlight the close similarity between certain quantities. For all the quantities defined above, the overlines may be ignored.  To obtain the formula for $\alpha$ from the formula for $\beta$,  replace $\overline{2}$ by $1$. 
To obtain the formula for $\delta$ from the formula for $\gamma$, replace the $\overline{3}$ by $4$ and change $(-3)^{i-1}$ to $(-3)^i$.
To obtain the formula for $\epsilon$ from the formula for $\zeta$, and also to obtain the formula for $\iota$ from the formula for $\kappa$,  replace each of the $\overline{2}$'s by $1$.
To obtain the formula for $\theta$ from the formula for $\eta$, and also to obtain the formula for $\mu$ from the formula for $\lambda$, replace each of the $\overline{3}$'s by $4$ and change $(-3)^{i-1}$ into $(-3)^i$.}

\end{theo}

\section{Proof of Theorem 1}

Let $P$ be a simple-2-column polyomino and let $P$ have at least two one-component columns. Then we define the \textit{lower pivot cell} of $P$ to be the cell which is the lower right neighbour of the bottom cell of the second-last (\textit{i.e.}, second-rightmost) among the one-component columns of $P$. We also define the \textit{upper pivot cell} of $P$ to be the upper right neighbour of the top cell of the second-last among the one-component columns of $P$.

Let $P$ be a simple-2-column polyomino and let $P$ have at least one two-component column. Then we define the \textit{lower inner pivot cell} of $P$ to be the upper right neighbour of the top cell of the lower component of the last among the two-component columns of $P$. We also define the \textit{upper inner pivot cell} of $P$ to be the lower right neighbour of the bottom cell of the upper component of the last among the two-component columns of $P$. 

Observe that the lower pivot cell of a simple-2-column polyomino $P$ is not necessarily contained in $P$. The same holds for the upper pivot cell, the lower inner pivot cell and the upper inner pivot cell of $P$. 

Let $\mathcal{S}$ denote the set of those simple-2-column polyominoes whose last (\textit{i.e.}, rightmost) column is a one-component column. It is convenient to first compute a generating function for the set $\mathcal{S}$, and thence a generating function for all simple-2-column polyominoes. So, let

\begin{displaymath}
A(q,t,w)=\sum_{P \in \mathcal{S}} q^{area \ of \ P} \cdot t^{{{the \ height \ of \ the \atop last \ column \ of \ P} \atop }} \cdot 
w^{{{the \ number \ of \atop two-component \ columns \ of \ P} \atop}}.
\end{displaymath}

Next, we define three generating functions in two variables, $q$ and $w$: Let $A_1=A(q,1,w)$, $B_1=\left[\frac{\partial A(q,t,w)}{\partial t}\right]_{with \ t=1}$ and $C_1=\frac{1}{2} \cdot \left[\frac{\partial^2 A(q,t,w)}{\partial t^2}\right]_{with \ t=1}$. \vspace{1mm}

The generating functions $G(q,w)$ and $A(q,t,w)$ are related by

\begin{equation}
\label{ekg}
G(q,w)=A_1 + \frac{q^2 w}{(1-q)^2} \cdot C_1.
\end{equation}

Henceforth the notation $A(q,t,w)$ will be abbreviated as $A(t)$.

In order to obtain a functional equation for the generating function $A(t)$, we are going to suitably partition the set $\mathcal{S}$. The blocks of the partition will be denoted $\mathcal{S}_{\alpha},\ \mathcal{S}_{\beta},\ldots ,\ \mathcal{S}_{\mu}$, and the parts of $A(t)$ coming from these blocks will be denoted $A_{\alpha}(t),\ A_{\beta}(t),\ldots ,\ A_{\mu}(t)$, respectively.

First, we define $\mathcal{S}_\alpha$ to be the set of those $P \in \mathcal{S}$ which have no other one-component column than the last column. We have $A_\alpha(t)  =\frac{qt}{1-qt}+\frac{q^4t^2w}{(1-q)^2(1-qt)^3}$.

Let $\mathcal{S}_\beta$ be the set of those $P \in \mathcal{S} \setminus \mathcal{S}_{\alpha}$ which have the following two properties: the second-last column is a one-component column, and the last column contains the lower pivot cell of $P$. We have $A_\beta(t)=\frac{qt}{(1-qt)^2} \cdot A_1$.

Let $\mathcal{S}_\gamma$ be the set of those $P \in \mathcal{S} \setminus \mathcal{S}_{\alpha}$ which have the following two properties: the second-last column is a one-component column, and the last column does not contain the lower pivot cell of $P$. We have $A_\gamma(t)=\frac{qt}{1-qt} \cdot B_1$.

Thus, $\mathcal{S}_{\beta} \cup \mathcal{S}_{\gamma}$ is the set of those $P \in \mathcal{S} \setminus \mathcal{S}_{\alpha}$ whose second-last column is a one-component column.

Let $\mathcal{S}_\delta$ be the set of those $P \in \mathcal{S} \setminus \mathcal{S}_{\alpha}$ which have the following three properties: the second-last column is a two-component column and the third-last column is (necessarily) a one-component column, the lower component of the second-last column and the third-last column have no edge in common, the lower pivot cell of $P$ is contained in the upper component of the second-last column. We have $A_\delta(t)=\frac{q^4t^2w}{(1-q)^3(1-qt)^3} \cdot A_1$.

Let $\mathcal{S}_\epsilon$ be the set of $P \in \mathcal{S} \setminus \mathcal{S}_{\alpha}$ having the following three properties: the second-last column is a two-component column and the third-last column is a one-component column, the lower component of the second-last column and the third-last column have no edge in common, and the lower pivot cell of $P$ is contained in the hole of the second-last column. We have $A_\epsilon(t) = \frac{q^4t^2w}{(1-q)^2(1-qt)^4} \cdot A_1 - \frac{q^4t^2w}{(1-q)^2(1-qt)^4} \cdot A(qt)$.

The definition of $\mathcal{S}_\zeta$ is obtained from the definition of $\mathcal{S}_\delta$ by writing the word ``upper" where the definition of $\mathcal{S}_\delta$ says ``lower", and by writing the word ``lower" where the definition of $\mathcal{S}_\delta$ says ``upper". We have $A_\zeta(t)=A_\delta(t)$.

The definition of $\mathcal{S}_\eta$ is obtained when the changes just described are made to the definition of $\mathcal{S}_\epsilon$ (instead of to the definition of $\mathcal{S}_\delta$). We have $A_\eta(t)=A_\epsilon(t)$.

Thus, $\mathcal{S}_{\delta} \cup \mathcal{S}_{\epsilon} \cup \mathcal{S}_{\zeta} \cup \mathcal{S}_{\eta}$ is the set of those $P \in \mathcal{S} \setminus \mathcal{S}_{\alpha}$ which, in addition to having a two-component second-last column, also have the property that $P \setminus (the \ last \ column \ of \ P)$ is not a polyomino. 

Let $\mathcal{S}_\theta$ be the set of $P \in \mathcal{S} \setminus \mathcal{S}_{\alpha}$ which have the following three properties: the second-last column is a two-component column and the third-last column is a one-component column, each of the two components of the second-last column has at least one edge in common with the third-last column, and the last column also has at least one edge in common with each of the two components of the second-last column. We have $A_{\theta}(t)= \frac{q^4t^2w}{(1-q)^2(1-qt)^3} \cdot B_1  - \frac{q^4t^2w}{(1-q)^2(1-qt)^4} \cdot A_1 + \frac{q^4t^2w}{(1-q)^2(1-qt)^4} \cdot A(qt)$.

Let $\mathcal{S}_\iota$ be the set of $P \in \mathcal{S} \setminus \mathcal{S}_{\alpha}$ which have the following four properties:

\begin{itemize}
\item the second-last column is a two-component column and the third-last column is a one-component column,
\item each of the two components of the second-last column has at least one edge in common with the third-last column,
\item the last column has at least one edge in common with the lower component of the second-last column, but does not have any edges in common with the upper component of the second-last column,
\item the last column does not contain the lower inner pivot cell of $P$.
\end{itemize}

We have $A_\iota(t) = \frac{q^3tw}{(1-q)^3(1-qt)} \cdot C_1$.

Let $\mathcal{S}_\kappa$ be the set of $P \in \mathcal{S} \setminus \mathcal{S}_{\alpha}$ having the following four properties:

\begin{itemize}
\item the second-last column is a two-component column and the third-last column is a one-component column,
\item each of the two components of the second-last column has at least one edge in common with the third-last column,
\item the last column has at least one edge in common with the lower component of the second-last column, but does not have any edges in common with the upper component of the second-last column,
\item the last column contains the lower inner pivot cell of $P$.
\end{itemize}

We have

\begin{eqnarray*}
A_\kappa(t) & = & \frac{q^3tw}{(1-q)^2(1-qt)^2} \cdot C_1 - \frac{q^4t^2w}{(1-q)^2(1-qt)^3} \cdot B_1 \\
& & \mbox{} + \frac{q^4t^2w}{(1-q)^2(1-qt)^4} \cdot A_1 - \frac{q^4t^2w}{(1-q)^2(1-qt)^4} \cdot A(qt) \ .
\end{eqnarray*}

The definition of $\mathcal{S}_\lambda$ is obtained from the definition of $\mathcal{S}_\iota$ by writing the word ``upper" where the definition of $\mathcal{S}_\iota$ says ``lower", and by writing the word ``lower" where the definition of $\mathcal{S}_\iota$ says ``upper". We have $A_\lambda(t)=A_\iota(t)$.

The definition of $\mathcal{S}_\mu$ is obtained when the changes just described are made to the definition of $\mathcal{S}_\kappa$ (instead of to the definition of $\mathcal{S}_\iota$). We have $A_\mu(t)=A_\kappa(t)$.

Thus, $\mathcal{S}_{\theta} \cup \mathcal{S}_{\iota} \cup \mathcal{S}_{\kappa} \cup \mathcal{S}_{\lambda} \cup \mathcal{S}_{\mu}$ is the set of those $P \in \mathcal{S} \setminus \mathcal{S}_{\alpha}$ which, in addition to having a two-component second-last column, also have the property that $P \setminus (the \ last \ column \ of \ P)$ is a polyomino. This means that $\mathcal{S}_{\delta} \cup \mathcal{S}_{\epsilon} \cup \ldots \cup \mathcal{S}_{\mu}$ is the set of all $P \in \mathcal{S} \setminus \mathcal{S}_{\alpha}$ whose second-last column is a two-component column. 

The sets $\mathcal{S}_\alpha, \ \mathcal{S}_\beta, \ldots, \ \mathcal{S}_\mu$ form a partition of the set $\mathcal{S}$. We have $A(t)=A_\alpha(t)+A_\beta(t)+\ldots+A_\mu(t)$, and the expressions for $A_\alpha(t), \ A_\beta(t), \ldots, \ A_\mu(t)$ are given above. Putting these things together, we get a functional equation for $A(t)$. It is convenient to write that functional equation as

\begin{eqnarray}
A(t) & = & \frac{qt}{1-qt} \cdot \left[ 1+B_1+\frac{2q^2w}{(1-q)^3} \cdot C_1 \right] 
+\frac{qt}{(1-qt)^2} \cdot \left[A_1+\frac{2q^2w}{(1-q)^2} \cdot C_1 \right]\nonumber \\
& & \mbox{} +\frac{q^4t^2w}{(1-q)^2(1-qt)^3} \cdot \left(1+\frac{2}{1-q} \cdot A_1 -B_1\right) 
+\frac{3q^4t^2w}{(1-q)^2(1-qt)^4} \cdot A_1 \nonumber \\
& & \mbox{} - \frac{3q^4t^2w}{(1-q)^2(1-qt)^4} \cdot A(qt). \label{eat}
\end{eqnarray}

We solved equation (\ref{eat}) by iteration, as is usually done in the upgraded Temperley method. The iteration ended in

\begin{eqnarray}
A(t) & = & \left\{\sum_{i=1}^{\infty} \frac{(-3)^{i-1}q^{i^2+2i-2}t^{2i-1}w^{i-1}}{(1-q)^{2i-2} \cdot \left[\prod_{k=1}^{i-1}(1-q^kt)\right]^4 \cdot (1-q^it)} \right\} \cdot \left[ 1+B_1+\frac{2q^2w}{(1-q)^3} \cdot C_1 \right] \nonumber \\
& & \mbox{} +\left\{\sum_{i=1}^{\infty} \frac{(-3)^{i-1}q^{i^2+2i-2}t^{2i-1}w^{i-1}}{(1-q)^{2i-2} \cdot \left[\prod_{k=1}^{i-1}(1-q^kt)\right]^4 \cdot (1-q^it)^2} \right\} \cdot \left[A_1+\frac{2q^2w}{(1-q)^2} \cdot C_1 \right] \nonumber \\
& & \mbox{} +\left\{\sum_{i=1}^{\infty} \frac{(-3)^{i-1}q^{i^2+3i}t^{2i}w^i}{(1-q)^{2i} \cdot \left[\prod_{k=1}^{i-1}(1-q^kt)\right]^4 \cdot (1-q^it)^3} \right\} \cdot \left(1+\frac{2}{1-q} \cdot A_1 -B_1\right) \nonumber \\
& & \mbox{} -\left\{\sum_{i=1}^{\infty} \frac{(-3)^i q^{i^2+3i}t^{2i}w^i}{(1-q)^{2i} \cdot \left[\prod_{k=1}^i(1-q^kt)\right]^4} \right\} \cdot A_1. \label{eatit}
\end{eqnarray}

From equation (\ref{eatit}), we got a system of three linear equations in three unknowns: $A_1$, $B_1$ and $C_1$. The first equation is just the case $t=1$ of (\ref{eatit}). The second equation is obtained by differentiating (\ref{eatit}) with respect to $t$ and then setting $t=1$. The third equation is obtained by differentiating (\ref{eatit}) twice with respect to $t$, and then setting $t=1$.

Once the linear system is solved, relation (\ref{ekg}) tells us how to obtain the sought-after generating function $G$. (To solve the linear system, we made use of the computer algebra package \textit{Maple}.)

\subsection{A corollary to Theorem 1}

By setting $w=0$, from Theorem 1 we obtain the well known result, discovered by Klarner \cite{Klarner}:

\begin{coro}
The area generating function for column-convex polygons with hexagonal cells is given by
\begin{displaymath}
G(q,0)=\frac{q(1-q)^3}{1-6q+10q^2-7q^3+q^4} \ .
\end{displaymath}
\end{coro}

\section{A bit of asymptotic analysis}

We write $[q^{n}] f$ to denote the coefficient of $q^{n}$ in a power series $f=f(q)$.

Formula (\ref{flag}) is very complicated, but still allows us to quickly compute Taylor polynomials of any reasonable degree. We actually computed the Taylor polynomial of degree $320$ for the function $G(q,1)$. It turned out that

\begin{eqnarray*}
G(q,1) & = & q+3q^2+11q^3+44q^4+186q^5+806q^6+3518q^7+15349q^8 \\
& & \mbox{}+66797q^9+289960q^{10}+1256274q^{11}+5435860q^{12}+\ldots \ .
\end{eqnarray*}

Then, for $n=2,\ 3,\ldots,\ 320$, we divided the coefficient $[q^n] G(q,1)$ by the coefficient $[q^{n-1}] G(q,1)$. This quotient gradually stabilizes, so that for $n=61,\ 62,\ldots,\ 320$ we have 

\begin{displaymath}
\left\{\frac{[q^n] G(q,1)}{[q^{n-1}] G(q,1)}\right\}_{rounded \ to \ 12 \ decimal \ places}=4.322382971063 \ .
\end{displaymath}

Simple concatenation arguments, first used by Klarner \cite{Klarner}, enable one to prove that $\mu = \lim_{n \to \infty}[q^n]^{1/n} = \sup_n [q^n]^{1/n}$. (The letter $\mu$ denotes the growth constant.) In this way, making use of $[q^{320}]$, we find that $\mu > 4.294676$, which is less than 1\% below the best numerical estimate. Unfortunately, finding upper bounds is much more difficult.

Next, for $n=1,\ 2,\ldots,\ 320$, we divided $[q^n] G(q,1)$ by $4.3223829710631654554^n$. Once again, the quotient stabilizes: for $n=54,\ 55,\ldots,\ 320$ we have $\frac{[q^n] G(q,1)}{4.3223829710631654554^n} \approx 0.127739087206 \:$. So, there is plenty of evidence that the coefficient $[q^n] G(q,1)$ (\textit{i.e.}, the number of $n$-celled simple-2-column polyominoes) has the asymptotic behaviour $[q^n] G(q,1) \sim 0.127739087206 \cdot 4.322382971063^n$. The dominant singularity of $G(q,1)$ is a simple pole, located at $0.231353863527 \:$.

For comparison, the area generating function for column-convex polyominoes is 

\begin{eqnarray*}
G(q,0) & = & q+3q^2+11q^3+42q^4+162q^5+626q^6+2419q^7+9346q^8 \\
& & \mbox{}+36106q^9+139483q^{10}+538841q^{11}+2081612q^{12}+\ldots \ .
\end{eqnarray*}

The number of $n$-celled column-convex polyominoes has the asymptotic behaviour\linebreak $[q^n] G(q,0) \sim 0.188419883819 \cdot 3.863130743243^n$. The dominant singularity of $G(q,0)$ is a simple pole, located at $0.258857405163 \:$.

As stated in \cite{Voege}, the area generating function for all polyominoes is

\begin{eqnarray*}
& & q+3q^2+11q^3+44q^4+186q^5+814q^6+3652q^7+16689q^8 \\
& & \mbox{}+77359q^9+362671q^{10}+1716033q^{11}+8182213q^{12}+\ldots \ 
\end{eqnarray*}

\noindent and the growth constant of all polyominoes is $5.183147\ldots\: $. Admittedly, the growth constant of all polyominoes is considerably greater than $4.322382\ldots \:$, the growth constant of simple-2-column polyominoes.

\section{Conclusion}

We have defined a class of polyominoes that interpolates between column-convex polyominoes and all polyominoes. The former have been solved, while the latter remain unsolved. The area generating function of column-convex polyominoes is just $\frac{q(1-q)^3}{1-6q+10q^2-7q^3+q^4}$. On the contrary, the area generating function of our interpolating model (called simple-2-column polyominoes) is a cumbersome rational function of $q$-series. The area generating functions of column-convex polyominoes and of simple-2-column polyominoes both have a simple pole singularity, located at $0.258857\ldots$ and $0.231353\ldots \:$, respectively. For all polyominoes, the corresponding singularity is at $q=q_c(\mathrm{polyomino})=0.192932\ldots \:$, and the singularity is of the form $const. \times \vert \mathrm{log}(q_c-q) \vert$, rather than a simple pole \cite{Voege}. We have also given the rigorous lower bound $\mu > 4.294676$, where $\mu$ is the growth constant of simple-2-column polyominoes.

In terms of possible extensions of this work, it is probably possible to compute the area generating function of \textit{simple-2-column}$_2$ \textit{polyominoes}. 
Here, by a simple-2-column$_2$ polyomino we mean a 2-column polyomino in which runs of two consecutive two-component columns are allowed, but it is forbidden for three consecutive columns to each have two connected components. One reason for doing this is that the growth factor $\mu$ is expected to be greater than that for simple-2-column polyominoes, and may set the benchmark in this regard. At present the situation is that for column-convex polygons the growth constant is $\mu=3.863130\ldots \:$, for simple-2-column polyominoes the growth constant is $\mu=4.322382\ldots \:$, and for all polyominoes the growth constant is $\mu=5.183147\ldots \:$. The polyomino model with a growth constant closest to the actual value for polyominoes is a directed model called {\it multi-directed polyominoes} \cite{BMR02} with a growth constant of $\mu \approx 4.587894\: $. It would be interesting to compute the area of simple-2-column$_2$ polyominoes to see if they had a growth constant closer still to that for polyominoes.

\section{Acknowledgement}

Svjetlan Fereti\'{c} wishes to thank Professor Anthony J. Guttmann for providing useful information, with appropriate references included, concerning the role that polyominoes play in theoretical physics.

\end{document}